\documentclass[12pt,twoside]{amsart}
\usepackage{amssymb}
\usepackage{latexsym}
\usepackage{amsfonts}
\usepackage{amsmath}
\usepackage{amssymb}
\usepackage{portland}

\newcommand{\K}{\mathcal K}

\newcommand{\E}{\mathcal E}
\newcommand{\cd}[1]{{\sf CD}(#1)}

\newtheorem{theorem}{Theorem}[section]
\newtheorem{proposition}[theorem]{Proposition}
\newtheorem{lemma}[theorem]{Lemma}

\theoremstyle{definition}

\newtheorem{definition}[theorem]{Definition}

\newcommand{\sy}[1]{{\sf S}_{#1}}

\newcommand{\sym}[1]{{\sf Sym}\,#1}
\renewcommand{\wr}{\,{\sf wr}\,}

\renewcommand{\leq}{\leqslant}
\renewcommand{\geq}{\geqslant}

\newcommand{\Z}{\mathbb Z}
\newcommand{\gl}[2]{{\sf GL}_{#1}(#2)}
\renewcommand{\sp}[2]{{\sf Sp}_{#1}(#2)}
\newcommand{\mat}[1]{{\sf M}_{#1}}

\newcommand{\aut}[1]{{\sf Aut}(#1)}
\newcommand{\out}[1]{{\sf Out}(#1)}

\newcommand{\cent}[2]{{\sf C}_{#1}(#2)}
\newcommand{\norm}[2]{{\sf N}_{#1}\left(#2\right)}

\newcommand{\nor}[1]{\mathbb N_{\F_{q^4}/\F_{q^2}}(#1)}

\newcommand{\psl}[2]{{\sf PSL}_{#1}(#2)}
\newcommand{\pomegap}[2]{{\sf P}\Omega^+_{#1}(#2)}

\newcommand{\alt}[1]{{\sf A}_{#1}}
\newcommand{\Alt}[1]{{\sf Alt}\,#1}

\newcommand{\F}{\mathbb F}

\begin{document}

\title{Transitive
simple subgroups of wreath products in product action}
\author{Robert W. Baddeley, Cheryl E. Praeger and Csaba Schneider}
\address[Baddeley]{32 Arbury Road,
  Cambridge CB4 2JE, UK}
\address[Praeger \& Schneider]{School of Mathematics and Statistics\\
The University of Western Australia\\
35 Stirling Highway 6009 Crawley\\
Western Australia}
\email{robert.baddeley@ntlworld.com, praeger@maths.uwa.edu.au,\newline
csaba@maths.uwa.edu.au\protect{\newline} {\it WWW:}
www.maths.uwa.edu.au/$\sim$praeger, www.maths.uwa.edu.au/$\sim$csaba}

\begin{abstract} A transitive simple subgroup of a finite symmetric group is
very rarely contained in a full wreath product in product action. All
such simple permutation groups are determined in this paper.  This remarkable
conclusion is reached after a definition and detailed examination of `Cartesian
decompositions' of the permuted set, relating them to certain `Cartesian systems
of subgroups'.  These concepts, and the bijective connections between them, are
explored in greater generality, with specific future applications in mind.
\end{abstract}

\date{12 January 2003}

\maketitle

\section{Introduction}\label{intro}

The main result of this paper is that a transitive simple subgroup of a finite
symmetric group is very rarely contained in a full
wreath product
in product action, so rarely that all such cases can be explicitly tabulated
here.  In other words, apart from a short list of exceptions, a simple subgroup
of a finite wreath product in product action can never be transitive.  A brief
summary of the product action of wreath products is provided at the beginning of
Section~\ref{sec2}.

\begin{theorem}\label{main}
Let $\Omega$ be a finite set, 
let $T<W<\sym\Omega$ such that $T$ is a finite
simple group, and $W$ is permutationally isomorphic to a wreath product
$\sym\Gamma\wr\sy\ell$ in product action. Then either $T$ is intransitive 
or $T$, $W$, and $|\Omega|$ are as in
Table~$\ref{maintable}$. Moreover, if $T$ is transitive, then $\norm{\sym\Omega}T$
is an almost simple group. 
\end{theorem}

\tiny
\begin{center}
\begin{table}[ht]
$$
\begin{array}{|c|c|c|c|}
\hline
 & T & W & |\Omega|\\
\hline\hline
1 &\alt 6& \sy 6\wr\sy 2& 36\\
\hline
2 & \mat{12}& \sy {12}\wr\sy 2 & 144 \\
\hline
3 & \pomegap 8q & \sy{(d/2)q^3(q^4-1)}\wr\sy 2
 & (d^2/4)q^6(q^4-1)^2\\
 && &d=(4,q^4-1) \\
\hline
4 &  \sp 4q,\ q\geq 4, q\, \mbox{ even}
& \sy{q^2(q^2-1)}\wr\sy 2  & q^4(q^2-1)^2\\
\hline
\end{array}
$$
\caption{ Transitive simple subgroups of wreath products}
\label{maintable}
\end{table}
\end{center}
\normalsize

This classification is reached after observing that, in Theorem~\ref{main}, the
set $\Omega$ can be identified with the Cartesian product $\Gamma^\ell$ such
that the action of $W$ is compatible with this identification.  In order to make
this idea precise, we introduce the concept of a `Cartesian decomposition' of a
set, and we also notice that $W$ can be viewed as the full stabiliser in
$\sym\Omega$ of a Cartesian decomposition of $\Omega$.  Hence we reduce the
problem of classifying the pairs $T,\ W$ in Theorem~\ref{main} to the problem of
classifying all Cartesian decompositions of finite sets that are invariant under
the action of a transitive, simple group of permutations.

Let $T$ be a finite simple group acting on a set $\Omega$. In the
classification of $T$-invariant Cartesian decompositions of $\Omega$
we use ideas that are familiar
from the elementary theory of permutation groups. Namely, we
investigate how the subgroup lattice of $T$ might reflect the existence of
a $T$-invariant Cartesian decomposition of $\Omega$. 
In Definition~\ref{csdef} we
define the concept of a `Cartesian system of subgroups', and in
Theorem~\ref{bijmain} we establish a one-to-one correspondence between the set
of $T$-invariant Cartesian decompositions of $\Omega$ and the set of
Cartesian systems with respect to a fixed element of $\Omega$.

The concepts of Cartesian decompositions and Cartesian systems, and the bijective 
connections between them, are explored in greater generality in
Sections~\ref{sec2}--\ref{bijsec}. Our motivation in doing so is to
provide with a theoretical background for a future investigation of
Cartesian decompositions that are invariant under a transitive
permutation group.

Some of the concepts we use may be new to most of our
readers. We define a permutation group
to be {\em innately transitive} if it has a transitive
minimal normal subgroup, and a transitive minimal normal subgroup of
an innately transitive group is referred to as a {\em plinth}. Most of
the results of this paper are expressed in the context of innately
transitive groups. 
The structure of innately transitive
groups is investigated in~\cite{bampra}. 
The problem of finding innately transitive subgroups of wreath
products in product action is studied more extensively in~\cite{bps}. 
Theorem~\ref{main} is  equivalent to the following result, which is
formulated in terms of innately transitive groups.
Here, a permutation group is \emph{quasiprimitive} if all of its minimal normal
subgroups are transitive.

\begin{theorem}\label{mainver}
Let $\Omega$ be a finite set, let $G<W<\sym\Omega$ such that  
$G$ is an innately transitive group with a
simple plinth $T$, and $W$ is permutationally isomorphic to a 
wreath product $\sym\Gamma\wr\sy\ell$ in product action. 
Then $T$ and $W$ are as in Table~$\ref{maintable}$, and 
$G$ is an almost simple quasiprimitive group.
\end{theorem}

Theorems~\ref{main} and \ref{mainver} are easy consequences of
Theorem~\ref{simth} as explained at the end of Section~\ref{simsec}.

\label{cdpage}A {\em Cartesian decomposition} of a finite set $\Omega$ is a
collection $\E$
of partitions $\Gamma_1,\ldots,\Gamma_\ell$ of $\Omega$ such that 
$$
|\gamma_1\cap\cdots\cap\gamma_\ell|=1\quad\mbox{for
 all}\quad\gamma_1\in\Gamma_1,\ldots,\gamma_\ell\in\Gamma_\ell.
$$
A Cartesian decomposition is said to be {\em homogeneous} if its elements
have the same size and this common size is at least~2. The number of partitions in a Cartesian
decomposition is called the {\em index}. A Cartesian decomposition is
said to be {\em non-trivial} if it has index at least~2. In this
paper, Cartesian decompositions are assumed to be non-trivial, unless
it is explicitly stated otherwise.
If $\E$ is a Cartesian decomposition of $\Omega$, then $\Omega$ can be
identified with the Cartesian product
$\prod_{\Gamma\in\E}\Gamma$. More information on Cartesian
decompositions is provided in~\cite{kov}, where a Cartesian decomposition $\E$
stabilised by a permutation group $G$ such that the elements of $\E$ form a single
$G$-orbit is said to be a system of product imprimitivity for $G$. 
A maximal
subgroup $W$ of $\sym\Omega$ or $\Alt\Omega$ 
is said to be {\em of product action type}, or
simply {\em PA type}, if $W$ is the full stabiliser of a non-trivial,
homogeneous Cartesian decomposition of $\Omega$. If 
a permutation group $G$ is contained in such a $W$, then we also say that $W$
is a {\em maximal overgroup of $G$ with product action type}, or
simply {\em PA type}. 

It is, in general, a difficult problem to describe maximal overgroups
with PA type of 
a transitive permutation group. In the case where $G$ itself is
primitive, this question is answered by~\cite{prae:inc},
but~\cite{bad:quasi} leaves this problem open for a quasiprimitive
$G$. Clearly our Theorem~\ref{mainver} gives a full classification of
the maximal overgroups of product action type for an innately
transitive permutation group $G$ with a simple plinth. This is
achieved by listing all non-trivial, homogeneous Cartesian decompositions stabilised by $G$. We
found that such decompositions
can be identified by information about the subgroups of the plinth. 
This motivates
the following definition.

\begin{definition}\label{csdef}
Let $M$ be a transitive 
permutation group on a set $\Omega$ and $\omega\in\Omega$. 
We say that a set $\{K_1,\ldots,K_\ell\}$ of subgroups
of $M$ is a {\em Cartesian system of subgroups} of $M$ with respect
to $\omega$ if 
\begin{eqnarray}
\label{csdef1}\bigcap_{i=1}^\ell K_i&=&M_\omega\quad\mbox{and}\\
\label{csdef2}K_i\left(\bigcap_{j\neq
i}K_j\right)&=&M\quad\mbox{for all}\quad i\in\{1,\ldots,\ell\}.
\end{eqnarray}
A Cartesian system is said to be {\em homogeneous} if its elements
are proper subgroups and they have the same size. A Cartesian system
is {\em non-trivial} if it has at least two subgroups. 
If $M$ is an abstract group and
$\K=\{K_1,\ldots,K_\ell\}$ is a set of subgroups
satisfying~\eqref{csdef2}, then $\K$ is called a Cartesian system of $M$.
\end{definition}

In this paper Cartesian systems are assumed to be non-trivial unless
it is explicitly stated otherwise.

For a permutation group $G\leq\sym\Omega$, 
let $\cd G$ denote the set of $G$-invariant
Cartesian decompositions of $\Omega$. 
Cartesian
systems provide a way of identifying the set $\cd G$ from
information internal to $G$. 

\begin{theorem}\label{bijmain}
Let $G$ be an innately transitive permutation group on $\Omega$ with plinth
$M$. Then for a fixed $\omega\in\Omega$ there is a one-to-one
correspondence between the set $\cd G$ and the set of
$G_\omega$-invariant Cartesian systems of $M$ with respect to
$\omega$.
\end{theorem} 

Theorem~\ref{bijmain} is an immediate consequence of
Theorem~\ref{bij} where an explicit one-to-one correspondence is
constructed.

The major results of this paper are
presented in Section~\ref{simsec}. There we study innately transitive
permutation groups with a non-abelian, simple plinth that preserve a
Cartesian decomposition of the underlying set. The main result of Section~\ref{simsec} gives rise to a complete description of maximal
overgroups  with product action type for such an innately
transitive group. Theorems~\ref{main} and \ref{mainver} follow 
immediately from Theorem~\ref{simth}(i),
where we give a detailed description of $G$-invariant 
homogeneous Cartesian decompositions of $\Omega$ for innately transitive groups
$G$ with a simple plinth $T$. In particular, Table~\ref{table}
contains the possibilities for $G$, $T$, $W$, $|\Omega|$, and the
isomorphism types of the subgroups in the associated Cartesian system,
as given by Theorem~\ref{bijmain}. Part~(ii) of Theorem~\ref{simth} gives a
detailed description of Cartesian decompositions $\E$ of
$\Omega$ with index at least~3 
that are invariant under the action of an innately transitive
group with a non-abelian, simple plinth. In Table~\ref{tablemult}, we
list the possibilities for the plinth, $|\Omega|$, the full stabiliser of
$\E$ in $\sym\Omega$, and the isomorphism types of the elements in the
corresponding Cartesian system. In the case where $G$ is
primitive, Theorem~\ref{simth} reduces
to~\cite[Proposition~6.1(ii)]{prae:inc}.
Problems similar to ours were also addressed
in~\cite{baum}.

Our notation concerning actions and 
permutation groups is standard. If $G$ is a
group acting on $\Omega$ and $\Delta$ is a subset of $\Omega$, then
$G_\Delta$ and $G_{(\Delta)}$ denote the setwise and the pointwise
stabilisers of $\Delta$, respectively. If $G_\Delta=G$ then $G^\Delta$
denotes the subgroup of $\sym\Delta$ induced by $G$. If
$\omega\in\Omega$, then $\omega^G$ denotes the $G$-orbit $\{\omega^g\
|\ g\in G\}$. 

\section{Cartesian decompositions}\label{sec2}

Let $\Gamma$ be a finite set, $L\leq\sym\Gamma$, $\ell\geq 2$ an
integer, and $H\leq\sy\ell$. The {\em
wreath product} $L\wr H$ is the semidirect product
$L^\ell\rtimes H$, where, for
$(x_1,\ldots,x_\ell)\in L^\ell$ and $\sigma\in\sy\ell$,
$(x_1,\ldots,x_\ell)^{\sigma^{-1}}=(x_{1^{\sigma}},\ldots,x_{\ell^{\sigma}})$. The product
action of $L\wr H$ is the action of $L\wr H$ on $\Gamma^\ell$ defined
by
$$
(\gamma_1,\ldots,\gamma_\ell)^{(x_1,\ldots,x_\ell)}=\left(\gamma_1^{x_1},\ldots,\gamma_\ell^{x_\ell}\right)\quad\mbox{and}\quad (\gamma_1,\ldots,\gamma_\ell)^{\sigma^{-1}}=(\gamma_{1^\sigma},\ldots,\gamma_{\ell^\sigma})
$$
for all $(\gamma_1,\ldots,\gamma_\ell)\in\Gamma^\ell$,
and $x_1,\ldots,x_\ell\in L$ and $\sigma\in H$.
The important properties of wreath products can be found in most
textbooks on permutation group theory, see for instance Dixon
and Mortimer~\cite{dm}.

The
full stabiliser $W$ in $\sym\Omega$ of a homogeneous Cartesian
decomposition $\E$ of $\Omega$ is
isomorphic to $\sym\Gamma\wr\sy\ell$
acting in product action on $\Gamma^\ell$ for $\Gamma\in\E$. Moreover,
if
$|\Gamma|\geq 3$ then $W$ is primitive on $\Omega$, and if
$|\Gamma|\geq 5$ then $W$ is a maximal subgroup of
$\sym\Omega$ or $\Alt\Omega$. 
As mentioned in Section~\ref{intro}, such maximal subgroups are usually referred to as {\em
maximal subgroups of product action type}. They
form one of several classes of primitive maximal
subgroups of $\sym\Omega$ and $\Alt\Omega$, identified by the O'Nan--Scott
Theorem; see~\cite{lps:maxsub}. Thus an important part of classifying
the primitive maximal subgroups of $\sym\Omega$ or $\Alt\Omega$ containing a given
(innately transitive) subgroup $G$ is finding all homogeneous Cartesian
decompositions 
of $\Omega$ that are stabilised by $G$. Our first result is that the
plinth must leave invariant each partition in such a Cartesian decomposition.

\begin{proposition}\label{minv}
If $G$ is an innately transitive group on a set $\Omega$ with plinth $M$ and 
$\E\in\cd G$, then $M_{(\E)}=M$. 
\end{proposition}
\begin{proof}
We let $\Gamma\in\E$ and show that each element of the $G$-orbit
$\Gamma^G$ is stabilised by $M$. 
Suppose that $\{\Gamma_1,\ldots,\Gamma_m\}$ is the $G$-orbit in $\E$
containing $\Gamma\in \E$. Set 
$$
\Sigma=\{\gamma_1\cap\cdots\cap\gamma_m\ |\
\gamma_1\in\Gamma_1,\ldots,\gamma_m\in\Gamma_m\}
$$
and
$$
\bar\Gamma_i=\{\{\sigma\in\Sigma\ |\ \sigma\subseteq\gamma\}\ |\
\gamma\in\Gamma_i\}\quad\mbox{for}\quad i=1,\ldots,m.
$$
Then it is a routine calculation to check that $\Sigma$ is a $G$-invariant partition of
$\Omega$, and that $\{\bar\Gamma_1,\ldots,\bar\Gamma_m\}$ is a $G$-invariant
Cartesian decomposition of $\Sigma$. Moreover, $|\bar
\Gamma_i|=|\Gamma_i|$ for all $i$, and since
$\Gamma_1,\ldots,\Gamma_m$ form a $G$-orbit, $|\bar
\Gamma_i|=|\bar\Gamma_j|$ for all $i$ and $j$. It is also easy to see
that if $g\in G_{(\bar\Gamma_i)}$ then $g\in G_{(\Gamma_i)}$. 
Since  $G_{(\bar\Gamma_1,\ldots,\bar\Gamma_m)}$ is a normal
subgroup of $G$ and $M$ is a minimal normal subgroup of $G$, 
either $M\leq G_{(\bar\Gamma_1,\ldots,\bar\Gamma_m)}$ or $M\cap
G_{(\bar\Gamma_1,\ldots,\bar\Gamma_m)}=1$. Suppose that $M\cap
G_{(\bar\Gamma_1,\ldots,\bar\Gamma_m)}=1$, so $M$ acts
on the set $\{\bar \Gamma_1,\ldots,\bar\Gamma_m\}$ faithfully. Therefore
$M$ is isomorphic to a subgroup
of $\sy m$. 
Note that $|\Sigma|=|\bar \Gamma_1|^m$, and let $p$ be a prime dividing $|\bar\Gamma_1|$. Then $p^m$
divides $|\Sigma|$. Since $M$ is transitive
on $\Sigma$, $p^m\mid|M|$. However, $M$ is isomorphic to a 
subgroup of $\sy m$, and so $p^m$ divides $m!$, which is a contradiction
to~\cite[Lemma~4.2]{prae:inc}. Hence $M\leq
G_{(\bar\Gamma_1,\ldots,\bar\Gamma_m)}$, that is, each
$\bar\Gamma_i$ is stabilised by $M$, and so is each $\Gamma_i$. Thus
$M$ stabilises $\Gamma$, and, since
$\Gamma$ was chosen arbitrarily, this shows that every element of
$\E$ is stabilised by $M$.
\end{proof}

\begin{lemma}\label{cslemma}
Let $M$ be a transitive subgroup of $\sym\Omega$ and let $\E\in\cd M$
such that $M_{(\E)}=M$. Suppose that
$\E=\{\Gamma_1,\ldots,\Gamma_\ell\}$, let $\omega\in\Omega$ be a fixed
element, and for $i=1,\ldots,\ell$ let $\gamma_i\in\Gamma_i$ be such that
$\omega\in\gamma_i$. Set $\K_\omega(\E)=\{K_1,\ldots,K_\ell\}$ where
$K_i=M_{\gamma_i}$ for $i=1,\ldots,\ell$. 
Then $\K_\omega(\E)$ is a Cartesian system of subgroups of $M$ with respect
to $\omega$. Moreover, if $\omega^m=\omega'$ for some $m\in M$, then $\K_{\omega'}(\E)=\K_\omega(\E)^m$.
\end{lemma}
\begin{proof}
Let us prove
that $\bigcap_{i=1}^\ell K_i=M_\omega$. 
Since
the $\Gamma_i$ are $M$-invariant partitions of $\Omega$, the
stabiliser of a point stabilises the block in $\Gamma_i$ that contains this
point. Hence $M_\omega\leq K_i$ for all $i$, and so
$M_\omega\leq\bigcap_i K_i$.
Now suppose 
$x\in \bigcap_{i} K_i$. Then $x$ stabilises
$\gamma_1,\ldots,\gamma_\ell$. Since $\E$ is a Cartesian
decomposition, $\gamma_1\cap\cdots\cap\gamma_\ell=\{\omega\}$, and so
$x$ stabilises $\omega$. 
Thus $x\in M_\omega$, and so  $\bigcap_{i}K_i=M_\omega$.

Now we prove that~\eqref{csdef2} also holds. We
may suppose without loss of generality that
$i=1$. Let $x\in M$,
$\delta_1=\gamma_1^x,\ldots,\delta_\ell=\gamma_\ell^x$, and
$\{\xi\}=\delta_1\cap\cdots\cap\delta_\ell$. If $\{\zeta\}=\delta_1\cap\gamma_2\cap\cdots\cap\gamma_\ell$ then 
the transitivity of $M$ on $\Omega$ implies that there exists $z\in M$
with $\xi^z=\zeta$ and so $\delta_1^z=\delta_1$,
$\delta_2^z=\gamma_2,\ldots,\delta_\ell^z=\gamma_\ell$, whence $\gamma_j^{xz}=\gamma_j$ for  $j=2,\ldots,\ell$ and $\gamma_1^{xzx^{-1}}=\gamma_1$, that is  $xz\in \bigcap_{j=2}^\ell K_j$ and $xzx^{-1}\in K_1$.
It follows that 
\[
x=(xzx^{-1})^{-1}(xzx^{-1}x)\in K_1 \left(\bigcap_{j=2}^\ell K_j\right),
\]
and we deduce that the first factorisation of~\eqref{csdef2} holds. The
other factorisations can be proved identically. Thus $\K_\omega(\E)$
is a Cartesian system of $M$ with respect to $\omega$.

If $m\in M$ and $\omega'=\omega^m$ then
$\{\omega'\}=\gamma_1^m\cap\cdots\cap\gamma_\ell^m$ and
$M_{\gamma_i^m}=M_{\gamma_i}^m$, which proves 
that $\K_{\omega'}(\E)=\K_\omega(\E)^m$. 
\end{proof}

If $M\leq\sym\Omega$ and $\E\in\cd M$ such that $M_{(\E)}=M$, then,
for a fixed $\omega\in\Omega$, we define the Cartesian system
$\K_\omega(\E)$  with respect to $\omega$
as in Lemma~\ref{cslemma}. The last result of this section establishes one direction of the
one-to-one correspondence in Theorem~\ref{bijmain}.

\begin{lemma}\label{eqactions}
Let $G$ be an innately transitive group with plinth $M$ acting on
$\Omega$, and let
$\omega\in\Omega$. If $\E\in\cd G$, then $M_{(\E)}=M$. Assume that 
$\K_\omega(\E)$ is the Cartesian
system of $M$ with respect to $\omega$. Then
$\K_\omega(\E)$ is invariant under conjugation by $G_\omega$, and the
$G_\omega$-actions on $\K_\omega(\E)$ and on $\E$ are equivalent.
\end{lemma}
\begin{proof} 
It follows from Proposition~\ref{minv} that $M_{(\E)}=M$, and so 
we can use Lemma~\ref{cslemma} to construct
$\K_\omega(\E)$ for $\omega$. 
Suppose that
$\E=\{\Gamma_1,\ldots,\Gamma_\ell\}$, and let
$\K_\omega(\E)=\{K_1,\ldots,K_\ell\}$ such that $K_i=M_{\gamma_i}$
where $\gamma_i$ is the unique element of $\Gamma_i$ containing
$\omega$. If $\Gamma_i,\ \Gamma_j\in\E$ and $g\in G_\omega$ such that
$\Gamma_i^g=\Gamma_j$ then $\omega^g=\omega$, and so 
$\gamma_i^g=\gamma_j$.
Hence
$$
K_i^g=\left(M_{\gamma_i}\right)^g=M_{\gamma_i^g}=M_{\gamma_j}=K_j,
$$
and so 
$\K_\omega(\E)$ is invariant under conjugation by $G_\omega$. This
argument also shows that the $G_\omega$-actions on $\E$ and on
$\K_\omega(\E)$ are equivalent. 
\end{proof}

\section{Cartesian systems}

In this section we summarise the most important properties of
Cartesian systems of abstract groups. The following lemma is useful
when working with Cartesian systems. 
If $\{K_1,\ldots,K_\ell\}$ is a Cartesian system for a group $M$ and
$I\subseteq\{1,\ldots,\ell\}$ then let $K_I$ denote the subgroup
$K_I=\bigcap_{i\in I}K_i$. We use the convention that if $I=\emptyset$
then $\bigcap_{i\in I}K_i=M$ for any collection $\{K_i\}_i$ of
subgroups in $M$.

\begin{lemma}\label{cscor}
Let $\{K_1,\dots,K_\ell\}$
be a (possibly trivial) Cartesian system for an abstract group $M$, and let $I$,
$J$ be subsets of
$\{1,\dots,\ell\}$.

(a) \ If $x_1,\dots,x_\ell\in M$, then $\bigcap_{i\in I}
K_ix_i$ is a coset modulo $K_I$.

(b) \  $|M:K_I|=\prod_{i\in I}|M:K_i|$.

(c) \ $K_IK_J=K_{I\cap J}$.  \end{lemma} 

\begin{proof} If an intersection of (right) cosets is nonempty then it is a
(right) coset modulo the intersection of the relevant subgroups.  The statement
of (a) above, and the simple proof below, make use of this fact.  We prove the
lemma by induction on $\ell$.  Notice that there is nothing to prove if
$\ell=1$.  Our inductive hypothesis is that $\ell>1$ and the lemma holds for all
Cartesian systems for $M$ which consist of fewer than $\ell$ subgroups.  Thus
(a) and (b) only have to be proved for the case $I=\{1,\dots,\ell\}$.  Put
$L=\bigcap_{i>1}K_i$, and note that $\{K_1,L\}$ is also a Cartesian system for
$M$ (that is, $K_1L=M$).

We also know from the inductive hypothesis that $\bigcap_{i>1} K_ix_i$ is a
coset modulo $L$, so for (a) it is sufficient to show that $K_1x_1\cap Ly$ is
never empty.  In order to show this we choose $z\in L$ such that
$K_1z=K_1x_1y^{-1}$; this is possible, as $K_1L=M$.  Then $K_1zy=K_1x_1$, and so
$zy\in K_1x_1$, and also $zy\in Ly$.  Hence $zy\in K_1x_1\cap Ly$, and
consequently $K_1x_1\cap Ly$ is non-empty.

For~(b), it is enough to show that
$|M:K_I|=|M:K_1||M:L|$, but this follows from 
$$
|M|=|K_1L|=|K_1||L|/|K_1\cap L|=|K_1||L|/K_I|.
$$
For an easy proof of (c) we first observe that
$$
|K_IK_J|=|K_I||K_J|/|K_I\cap K_J|.
$$
It is obvious that  $K_IK_J\subseteq K_{I\cap J}$ and, as 
$K_I\cap K_J=K_{I\cup J}$, one can calculate from (b) and the last display
that
$|K_IK_J|=|K_{I\cap J}|$. This completes the proof of the lemma.
\end{proof}

Note that, in Lemma~\ref{cscor}(a), if we choose $x$ to be any element
of $\bigcap_{i\in I}K_ix_i$, then $K_ix_i=K_ix$ holds, for all $i\in I$.

\section{Cartesian systems and Cartesian decompositions}\label{bijsec}

In a transitive group $M\leq\sym\Omega$, a subgroup $K$
satisfying $M_\omega\leq K\leq M$ for some $\omega\in\Omega$
determines an $M$-invariant partition of $\Omega$ comprising the
$M$-translates of the $K$-orbit $\omega^K$. 

\begin{lemma}\label{cs2cd}
Let $G$ be an innately transitive group on $\Omega$ 
with plinth $M$, and let $\omega$ be a fixed element of
$\Omega$. Suppose that $\K=\{K_1,\ldots,K_\ell\}$ is a $G_{\omega}$-invariant 
Cartesian system of subgroups of $M$ with respect to $\omega$, and
let $\Gamma_1,\ldots,\Gamma_\ell$ be the $M$-invariant partitions of
$\Omega$
determined by $K_1,\ldots,K_\ell$, respectively. Then
$\E=\{\Gamma_1,\ldots,\Gamma_\ell\}$ is a $G$-invariant Cartesian
decomposition of $\Omega$, such that $\K_\omega(\E)=\K$. 
Moreover, if $M$ is non-abelian and the Cartesian system
$\{K_1,\ldots,K_\ell\}$ is homogeneous, then 
the stabiliser $W$ in $\sym\Omega$ of
$\E$ is a maximal subgroup of $\sym\Omega$ or $\Alt\Omega$ such that $G\leq W$.
\end{lemma}
\begin{proof}
As $M_\omega\leq K_i\leq M$, 
each $\Gamma_i$ is an $M$-invariant partition of $\Omega$. For
$i=1,\ldots,\ell$ let $\gamma_i$ be the unique element of $\Gamma_i$
containing $\omega$. In order to 
prove that $\E$ is a Cartesian decomposition, we only have to show that
$$
\left|\bigcap_{i=1}^\ell\delta_i\right|=1\quad\mbox{whenever}\quad
\delta_1\in\Gamma_1,\ldots,\delta_\ell\in\Gamma_\ell.
$$
To see this, choose 
$\delta_1\in\Gamma_1,\ldots,
\delta_\ell\in\Gamma_\ell$. Now
$\delta_i=\gamma_i^{x_i}$ for some $x_i\in M$, and by Lemma~\ref{cscor}(a),
there exists some $x\in M$ such that $K_ix_i=K_ix$ for $i=1,\ldots,\ell$. 
Then
\begin{multline*}
\delta_i=\gamma_i^{x_i}=\{\omega^k\,|\,k\in K_i\}^{x_i}=\{\omega^{k'}\,|\,k'\in K_ix_i\}\\
=\{\omega^{k'}\,|\,k'\in K_ix\}=\{\omega^{k}\,|\,k\in K_i\}^x=\gamma_i^x.
\end{multline*}
Thus 
$$
\bigcap_{i=1}^\ell\delta_i=\bigcap_{i=1}^\ell\gamma_i^x=\left(\bigcap_{i=1}^\ell\gamma_i\right)^x,
$$
and therefore we only have to prove that $|\bigcap_{i=1}^\ell\gamma_i|=1$. 
Note that $\omega\in\gamma_i$ for  $i=1,\ldots,\ell$. 
Suppose that $\omega'\in\gamma_1\cap\ldots\cap\gamma_\ell$ for some
$\omega'\in\Omega$. Then there is some $x\in M$ such that
$\omega^x=\omega'$. Then $x$ must stabilise
$\gamma_1,\ldots,\gamma_\ell$, and hence 
$x\in K_i$ for all $i=1,\ldots,\ell$. Since $\bigcap_{i=1}^\ell K_i=M_\omega$,
it follows that $x\in M_\omega$, and so $\omega^x=\omega$. Thus
$\bigcap_{i=1}^\ell\gamma_i=\{\omega\}$, and $\E$ is a Cartesian
decomposition. 

Since each $\Gamma_i$ is an $M$-invariant partition of $\Omega$, $\E$
is invariant under $M$. Since $\{K_1,\ldots,K_\ell\}$ is
$G_\omega$-invariant, $\E$ is also $G_\omega$-invariant, and so $\E$
is $MG_\omega$-invariant. Since $M$ is transitive,
$MG_\omega=G$. Therefore $\E$ is $G$-invariant. Note that 
$$
\K=\{M_{\gamma_1},\ldots,M_{\gamma_\ell}\}\quad\mbox{and}\quad
K_\omega(\E)=\{M_{\gamma_1},\ldots,M_{\gamma_\ell}\}.
$$
Thus $\K=\K_\omega(\E)$, as required.

Since $M$ is non-abelian, $M$ is a direct product of isomorphic
non-abelian, simple groups. Hence for $i=1,\ldots,\ell$, the group
$M^{\Gamma_i}$ is also isomorphic to a direct product of non-abelian
simple groups. Moreover, $M^{\Gamma_i}$ is transitive and faithful on
$\Gamma_i$, and so $|\Gamma_i|\geq 5$ for all $i$. As
$\{K_1,\ldots,K_\ell\}$ is homogeneous, $\E$ is also homogeneous and
$W$ is permutationally isomorphic to $\sym\Gamma\wr\sy\ell$ in product
action for some set $\Gamma$ and $\ell\geq 2$. Hence the results
of~\cite{lps:maxsub} show that $W$ is a maximal subgroup of
$\sym\Omega$ if $W\not\leq\Alt\Omega$, and $W$ is a maximal subgroup
of $\Alt\Omega$ otherwise. Since $\E$ is $G$-invariant, clearly $G\leq W$.
\end{proof}

\begin{theorem}\label{bij}
Let $G$ be an innately transitive group on $\Omega$ with plinth
$M$. For a fixed $\omega\in\Omega$ 
the map $\E\mapsto\K_\omega(\E)$ is a bijection between the set $\cd G$ 
and the set of $G_\omega$-invariant
Cartesian systems of subgroups of $M$ with respect to $\omega$.
\end{theorem}
\begin{proof}
Let $\mathcal C$ denote the set of $G_\omega$-invariant
Cartesian systems of subgroups of $M$ with respect to $\omega$.
In Lemma~\ref{cslemma}, we explicitly constructed a map
$\Psi: \cd G\rightarrow \mathcal C$ for which
$\Psi(\E)=\K_\omega(\E)$. We claim that $\Psi$ is a
bijection. Let $\K\in\mathcal C$,  let
$\Gamma_1,\ldots,\Gamma_\ell$ be the $M$-invariant partitions determined
by the elements $K_1,\ldots,K_\ell$ of $\K$, and let
$\E=\{\Gamma_1,\ldots,\Gamma_\ell\}$. We proved in Lemma~\ref{cs2cd} that
$\E$ is a $G$-invariant Cartesian decomposition of $\Omega$ such that
$\K_\omega(\E)=\K$. Hence $\Psi$ is surjective.

Suppose now that $\E_1,\ \E_2\in\cd G$ is such that
$\Psi(\E_1)=\Psi(\E_2)$ and let $\K$ denote this common Cartesian system. Let
$\E$ be the set of $M$-invariant
partitions determined by the elements of $\K$. Then, by the definition
of $\Psi(\E_i)$ in Lemma~\ref{cslemma}, $\E_1=\E$ and
$\E_2=\E$. Thus $\Psi$ is
injective, and so $\Psi$ is a bijection.
\end{proof}

Theorem~\ref{bijmain} is an immediate consequence of the previous
result.

\section{Some factorisations of finite simple groups}

To prove Theorem~\ref{main} 
we need first to prove some results about factorisations of certain finite
simple groups. If $G$ is a group and $A,\ B\leq G$ such that $G=AB$,
then we say that the expression $G=AB$ or the set $\{A,B\}$ 
is a {\em factorisation} of $G$. In~\cite{bad:fact} full
factorisations of almost simple groups were classified up to the
following equivalence relation. The factorisations $G=A_1B_1$ and
$G=A_2B_2$ of a group $G$ are said to be {\em equivalent} if there are
$\alpha\in\aut G$, and $x,\ y\in G$ such that
$\{A_1,B_1\}=\{A_2^{\alpha x},B_2^{\alpha y}\}$. The following lemma
shows that this equivalence relation can be expressed in a simpler way.

\begin{lemma}\label{equivlemma}
Let $G$ be a group.
\begin{enumerate}
\item[(i)] If $G=AB$ for some $A,\ B\leq G$, then the conjugation
action of $A$ is transitive on the conjugacy class $B^G$, and $B$ is
transitive on $A^G$. 
\item[(ii)] The factorisations $G=A_1B_1$ and
$G=A_2B_2$ of $G$ are equivalent if and only if there is
$\beta\in\aut G$ such that $\{A_1,B_1\}=\{A_2^\beta,B_2^\beta\}$. 
\end{enumerate}
\end{lemma}
\begin{proof}
(i) 
As $AB=G$, we also have $A\norm G{B}=G$. Since $\norm
GB$ is a point stabiliser for the conjugation action of $G$ on the
conjugacy class $B^G$, we obtain that $A$ is a transitive subgroup of
$G$ with respect to this action. Similar argument shows that $B$ is
transitive by conjugation on $A^G$. 

(ii) It is clear that if there is $\beta\in\aut G$ such that
$\{A_1,B_1\}=\{A_2^{\beta},B_2^{\beta}\}$ then the two
factorisations in the lemma are equivalent. 
Suppose that $G=A_1B_1$ and $G=A_2B_2$ are equivalent factorisations. 
By assumption, there is $\alpha\in\aut G$ and $x,\ y\in G$ such that
$\{A_1,B_1\}=\{A_2^{\alpha x},B_2^{\alpha y}\}$. 
Then we have $A_1^G=(A_2^{\alpha })^G$ and $B_1^G=(B_2^{\alpha
})^G$, or $A_1^G=(B_2^{\alpha
})^G$ and $B_1^G=(A_2^{\alpha
})^G$. Suppose without loss of generality that $A_1^G=(A_2^{\alpha
})^G$ and $B_1^G=(B_2^{\alpha
})^G$. 
Since $A_1$ and $A_2^\alpha$ are conjugate, there is some
$g\in G$ such that $A_1^g=A_2^\alpha$, and $B_1^g$ is conjugate to
$B_2^\alpha$. 
As $G=(A_1B_1)^g=A_1^gB_1^g$, we have that $A_1^g$ is transitive by
conjugation on $(B_1^g)^G=B_1^G$. Hence there is some $a\in A_1^g$ such that
$A_1^{ga}=A_1^g=A_2^\alpha$, and $B_1^{ga}=B_2^\alpha$. Hence
$A_1=A_2^{\alpha a^{-1}g^{-1}}$ and $B_1=B_2^{\alpha a^{-1}g^{-1}}$. Thus we
may take $\beta$ as $\alpha$ followed by the inner automorphism
corresponding to $a^{-1}g^{-1}$. 
\end{proof}

If $G$ is a group and $A$ and $B$ are subgroups then
let 
$$
\norm G{\{A,B\}}=\{g\in
G\ |\ \{A^g,B^g\}=\{A,B\}\}.
$$ 
In the proof of the following result we use the following simple fact,
called Dedekind's modular law. If $K,\ L,\ H$ are subgroups of a group
$G$ such that $K\leq L$, then
\begin{equation}\label{dedekind}
(HK)\cap L=(H\cap L)K.
\end{equation}

\begin{lemma}\label{factlem}
Let $T$ be a finite simple group and $A,\ B$ proper subgroups of $T$ such
that $|A|=|B|$ and $T=AB$. Then the following hold.
\begin{enumerate}
\item[(i)] The isomorphism types of 
$T$, $A$, and $B$ are as in
Table~$\ref{isomfact}$, and $A$, $B$ are maximal subgroups of $T$.
\item[(ii)] There is an automorphism $\vartheta\in\aut T$ such that
$\vartheta$ interchanges $A$ and $B$.
\item[(iii)] The group $A\cap B$
is self-normalising in $T$. 
\item[(iv)] If $T$ is as in row~$1$, $2$, or~$4$ of
Table~$\ref{isomfact}$, then 
$$
\norm {\aut T}{A\cap B}=\norm{\aut
T}{\{A,B\}}=N, 
$$
say, and moreover $TN=\aut T$.
\end{enumerate}
\end{lemma}
\begin{center}
\begin{table}[ht]
$$
\begin{array}{|l|c|c|}
\hline
 & T & A,\ B\ \\
\hline
1 & \alt 6&\alt 5\\
\hline
2 & \mat{12}&\mat{11}
 \\
\hline
3 & \pomegap 8q &\Omega_7(q)\\
\hline
4 & \sp 4q,\ q\geq 4\mbox{ even}
& \sp 2{q^2}.2  \\
\hline
\end{array}
$$
\caption{Factorisations of finite simple groups in Lemma~\ref{factlem}}
\label{isomfact}
\end{table}
\end{center}
\begin{proof}
(i) Note that, since $|A|=|B|$, the factorisation $T=AB$ is a full
factorisation of $T$, that is, the sets of primes dividing $|T|$,
$|A|$, and $|B|$ are the same. It was proved in~\cite{bad:fact}, that
$T$, $A$, and $B$ are as in~\cite[Table~I]{bad:fact}. It is easy to see that the only possibilities
where $|A|=|B|$ are those in Table~\ref{isomfact}, and it follows
that in these cases $A$ and $B$ are maximal subgroups of $T$.

(ii) In each line of Table~\ref{isomfact}, the groups $A$
and $B$ are not conjugate, but there is an
outer automorphism $\sigma\in\aut T$ which swaps the conjugacy classes $A^T$ and
$B^T$ (see the Atlas~\cite{atlas} for $T\cong \alt 6,\ \mat{12}$,
\cite{Kleidman} for $T\cong\pomegap 8q$, and \cite[page~155]{bad:fact} for $T\cong\sp 4q$). 
By Lemma~\ref{equivlemma}(i), the group $A$ is
transitive in its conjugation action on the conjugacy class $B^T$ 
and $B$ is transitive on $A^T$. Thus there is an
element $a\in A$ such that $A^{\sigma a}=B$ and $B^{\sigma a}$ is
conjugate in $T$ to $A$. Since $B$ is transitive on $A^T$, there
is an element $b\in B$ such that $A^{\sigma ab}=B$ and $B^{\sigma ab
}=A$.  Therefore we can take $\vartheta$ as $\sigma$ followed by the
inner automorphism induced by the element $ab$.

(iii) Set $C=A\cap B$. First we prove that $C$ is self-normalising in
$T$. If $T$ is isomorphic to $\alt 6$ or $\mat{12}$ then the
information given in the Atlas~\cite{atlas} shows that if $N$ is a proper subgroup
of $T$ properly containing $C$, then $N$ is isomorphic to $A$ or $B$. In all
cases $A$ and $B$ are simple, and so $\norm TC=C$. 
If $T\cong\pomegap 8q$ then we obtain from~\cite[3.1.1(vi)]{Kleidman} that $C\cong {\sf G}_2(q)$ and~\cite[3.1.1(iii)]{Kleidman} yields
that $\norm TC=C$. 

Now let $T\cong\sp 4q$ for $q\geq 4$, $q$ even. 
In this case  $A\cong B\cong
\sp 2{q^2}\cdot 2$. Consider the fields
$\F_q$ and $\F_{q^2}$ as subfields of the field $\F_{q^4}$ and 
consider the field
$\F_{q^4}$ as a 4-dimensional vector space $V$ over $\F_q$. Let
$\mathbb N_{\F_{q^4}/\F_{q^2}}:\F_{q^4}\rightarrow \F_{q^2}$ and
$\textrm{Tr}_{\F_{q^2}/\F_{q}}:\F_{q^2}\rightarrow\F_{q}$ denote the norm and the trace
map, respectively. For the basic properties of these maps
see~\cite[2.3]{ln}. 
Using the fact that $\nor x=x^{q^2+1}$ for all $x\in\F_{q^4}$, we
obtain that $x\mapsto \nor x$ is an $\F_{q^2}$-quadratic form on $V$,
such that $(x,y)\mapsto \nor{x+y}+\nor x+\nor y$ is a non-degenerate, symmetric,
$\F_{q^2}$-bilinear form with Witt defect~1 
(we recall that $q$ is a 2-power). Hence $Q=\textrm{Tr}_{\F_{q^2}/\F_{q}}\circ\mathbb N_{\F_{q^4}/\F_{q^2}}$
is an $\F_{q}$-quadratic form $V\rightarrow\F_{q}$, and
$f(x,y)=Q(x+y)+Q(x)+Q(y)$, is a non-degenerate, symmetric
$\F_q$-bilinear form on $V$ with Witt defect~1. Then without loss of generality we may assume that
$T$ is the stabiliser of $f$ in $\gl 4q$, $A$ consists of elements of $T$ that are
$\F_{q^2}$-semilinear, and $B$ is the stabiliser of $Q$.

For $a\in V\setminus\{0\}=\F_{q^4}^*$, define the map
$s_a:x\mapsto xa$. Then  it is well-known that $S=\{s_a\ |\
a\in\F_{q^4}^*\}$ is a cyclic subgroup of $\gl 4q$. A generator of $S$
is called a Singer cycle; see Satz~II.7.3 in
Huppert~\cite{huppert}. Let $Z$ denote the subgroup $\{s_a\ |\
\mathbb N_{\F_{q^4}/\F_{q^2}}(a)=1\}$ of $S$. Since the restriction of $\mathbb N_{\F_{q^4}/\F_{q^2}}$ to $\F_{q^4}^*$ is an
epimorphism $\mathbb
N_{\F_{q^4}/\F_{q^2}}:\F_{q^4}^*\rightarrow\F_{q^2}^*$, and $Z$ is
the
kernel of this epimorphism, we have that $|Z|=q^2+1$. If $\sigma$ is
the Frobenius automorphism $x\mapsto x^q$ of $\F_{q^4}$ then
$(s_a)^\sigma=s_{a^\sigma}$ 
for all $s_a\in S$. Therefore $\sigma$
normalises $S$, and, since $S$ is cyclic, $\sigma$ also normalises
$Z$. We claim that
$C=Z\left<\sigma\right>$. Since $T=AB$, $|C|=4(q^2+1)$, and hence 
it suffices to
prove that $Z\left<\sigma\right>\leq C$. It is clear that $\sigma$ is
$\F_{q^2}$-semilinear, and so $\sigma\in A$. Also
\begin{multline*}
Q(\sigma(x))=Q(x^q)=\textrm{Tr}_{\F_{q^2}/\F_{q}}\left(\mathbb N_{\F_{q^4}/\F_{q^2}}\left(x^q\right)\right)=\textrm{Tr}_{\F_{q^2}/\F_{q}}\left(\mathbb N_{\F_{q^4}/\F_{q^2}}\left(x\right)^q\right)\\=\textrm{Tr}_{\F_{q^2}/\F_{q}}\left(\mathbb N_{\F_{q^4}/\F_{q^2}}\left(x\right)\right)=Q(x).
\end{multline*}
Therefore $\sigma\in B$, and so $\sigma\in C$. Let $a\in\F_{q^4}$ such that
$\mathbb N_{\F_{q^4}/\F_{q^2}}(a)=1$. Then
\begin{multline*}
Q(s_a(x))=Q(xa)=\textrm{Tr}_{\F_{q^2}/\F_{q}}\left(\mathbb N_{\F_{q^4}/\F_{q^2}}(xa)\right)\\=\textrm{Tr}_{\F_{q^2}/\F_{q}}\left(\mathbb N_{\F_{q^4}/\F_{q^2}}(x)\mathbb N_{\F_{q^4}/\F_{q^2}}(a)\right)=\textrm{Tr}_{\F_{q^2}/\F_{q}}\left(\mathbb N_{\F_{q^4}/\F_{q^2}}(x)\right)=Q(x).
\end{multline*}
Thus $s_a\in B$. Since $s_a$ is also $\F_{q^2}$-linear, we obtain $s_a\in
C$. Hence $C=Z\left<\sigma\right>$. 

We will now prove that $C$ is self-normalising in $T$.
First notice that $q^2+1$ is divisible by an odd prime $r$ such that $r\nmid (q^2-1)$. Hence there is
a unique subgroup $R$ in $Z$ with order $r$. Since $Z$ is the
commutator subgroup of $C$, it is a characteristic subgroup of $C$.
Also $R$ is the unique
subgroup of $Z$ with order $r$, and so $R$ is characteristic in $Z$. Thus $R$ is
 characteristic in
$C$ and $\norm TC$ must
normalise $R$. 
By Satz~II.7.3 of Huppert~\cite{huppert}, $\norm{\gl 4q}R=SC=S\left<\sigma\right>$. 

Let us now determine how much of $S\left<\sigma\right>$ is contained in $T$. 
Since $\textrm{Tr}_{\F_{q^2}/\F_q}$ is additive,
\begin{multline*}
f(x^q,y^q)=\textrm{Tr}_{\F_{q^2}/\F_q}\left({\mathbb
N}_{\F_{q^4}/\F_{q^2}}(x^q+y^q)+{\mathbb
N}_{\F_{q^4}/\F_{q^2}}(x^q)+{\mathbb
N}_{\F_{q^4}/\F_{q^2}}(y^q)\right)\\=\textrm{Tr}_{\F_{q^2}/\F_q}\left(\left({\mathbb 
N}_{\F_{q^4}/\F_{q^2}}(x+y)+{\mathbb
N}_{\F_{q^4}/\F_{q^2}}(x)+{\mathbb
N}_{\F_{q^4}/\F_{q^2}}(y)\right)^q\right)=f(x,y),
\end{multline*}
and hence the cyclic subgroup $\left<\sigma\right>$ is in $T$. Using~\eqref{dedekind}, we have $(S\left<\sigma\right>)\cap T=(S\cap
T)\left<\sigma\right>$. Thus we need to compute $S\cap T$.
If $x\in \F_{q^4}^*$ such that $f(xa,xb)=f(a,b)$
for all $a,\ b\in V$ then 
\begin{multline}\label{normeqhere}
\textrm{Tr}_{\F_{q^2}/\F_q}\left(\mathbb N_{\F_{q^4}/\F_{q^2}}(a+b)+\mathbb
N_{\F_{q^4}/\F_{q^2}}(a)+\mathbb
N_{\F_{q^4}/\F_{q^2}}(b)\right)\\=\textrm{Tr}_{\F_{q^2}/\F_q}\left(\mathbb
N_{\F_{q^4}/\F_{q^2}}(xa+xb)+\mathbb
N_{\F_{q^4}/\F_{q^2}}(xa)+\mathbb N_{\F_{q^4}/\F_{q^2}}(xb)\right)\\=\textrm{Tr}_{\F_{q^2}/\F_q}\left(\mathbb
N_{\F_{q^4}/\F_{q^2}}(x)\left(\mathbb
N_{\F_{q^4}/\F_{q^2}}(a+b)+\mathbb N_{\F_{q^4}/\F_{q^2}}(a)+\mathbb N_{\F_{q^4}/\F_{q^2}}(b)\right)\right). 
\end{multline}
As observed above,
$$
(u,v)\mapsto\mathbb N_{\F_{q^4}/\F_{q^2}}(u+v)+\mathbb
N_{\F_{q^4}/\F_{q^2}}(u)+\mathbb N_{\F_{q^4}/\F_{q^2}}(v)
$$ 
is a non-degenerate, symmetric, $\F_{q^2}$-bilinear form, and so it maps
$V\times V$ onto $\F_{q^2}$. 
Hence~\eqref{normeqhere} shows that 
$y=\mathbb N_{\F_{q^4}/\F_{q^2}}(x)$ has the property that
$\textrm{Tr}_{\F_{q^2}/\F_q}(yu)=\textrm{Tr}_{\F_{q^2}/\F_q}(u)$ for
all $u\in \F_{q^2}$, that is, $yu+y^qu^q=u+u^q$, for all
$u\in\F_{q^2}$. Thus $(yu+u)^q=yu+u$. Hence $u(y+1)\in\F_q$ for all $u\in\F_{q^2}$,
and consequently $y=1$. Thus if the map $s_x$ preserves $f$ then
$\mathbb N_{\F_{q^4}/\F_{q^2}}(x)=1$. On the other hand from~\eqref{normeqhere} it is clear that if
$\mathbb N_{\F_{q^4}/\F_{q^2}}(x)=1$ then multiplication by $x$ preserves $f$. Since the norm is
a group epimorphism $\mathbb
N_{\F_{q^4}/\F_{q^2}}:\F_{q^4}^*\rightarrow\F_{q^2}^*$ it follows that
the elements of norm~1 form  a cyclic group of order $q^2+1$. Hence
$S\cap T=Z$ and $\norm {\gl 4q}C\cap T=C$, that is, $C$ is
self-normalising in $T$.

(iv) Finally we assume that $T$ is as in row~1, 2, or~4 of
Table~\ref{isomfact}, and we prove the assertion that $N_1=N_2$, where
$N_1=\norm
{\aut T}{C}$ and $N_2=\norm{\aut T}{\{A,B\}}$.
It is clear that
$N_2\leq N_1$, and so we only have to prove $|N_1|\leq |N_2|$. Since $A$ and $B$ are not
conjugate in $T$, we have that $N_2\cap T=\norm TA\cap\norm TB=A\cap
B=C$, and, since $C$ is self-normalising in $T$, we also have $N_1\cap
T=C$. Thus it suffices to prove that
$TN_1\leq TN_2$, which follows immediately once we show that
$TN_2=\aut T$.  Since $N_2$ interchanges $A$ and $B$, we have
that $N_2=(\norm{\aut T}A\cap\norm{\aut T}B)\left<\vartheta\right>$
where $\vartheta\in\aut T$ is as in~(ii). If $T\cong \alt 6$ then 
$|\norm {\aut T}A:\norm TA|=|\norm{\aut T}B:\norm TB|=2$, and so
$TN_2=\aut T$ (see~\cite{atlas}). If $T\cong\mat{12}$ then $\norm{\aut T}A=\norm TA$ and
$\norm {\aut T}B=\norm TB$, and so $TN_2=\aut T$ (see~\cite{atlas}).
If $T\cong\sp 4q$ then the field automorphism group $\Phi$ normalises
$A$ and $B$. If $\vartheta\in\aut T$ is as in~(ii), then $\aut
T=T\Phi\left<\vartheta\right>$, and so
we obtain that $TN_2=\aut T$. Hence if $T$ is as in row~1, 2, or~4 of
Table~\ref{isomfact}, then 
$TN_2=\aut T$, and $TN_1\leq TN_2$ clearly holds. Thus $N_1=N_2$ follows.
\end{proof}

We recall a couple of facts about automorphisms  of  
$\pomegap 8q$. Let $T=\pomegap 8q$. Then, as shown in~\cite[pp.~181--182]{Kleidman},
$\aut T=\Theta\rtimes\Phi$, where $\Phi$ is the group of  field
automorphisms of $T$, and $\Theta$ is a certain subgroup of $\aut T$
containing the commutator subgroup $\aut T'$. We also have $\out
T=\aut T/T=\Theta/T\times\Phi T/T$, and $\Theta/T\cong\sy m$ where
$m=3$ for even $q$, and $m=4$ for odd $q$. Let $\pi:\Theta\rightarrow\sy m$
denote the natural epimorphism. The following lemma derives the
information about $\pomegap 8q$ similar to that in Lemma~\ref{factlem}(iv).

\begin{lemma}\label{pomegalem}
Let $T=\pomegap 8q$, let $A$, $B$ be
subgroups of $T$ such that $A,\ B\cong\Omega_7(q)$ and
$AB=T$, and set $C=A\cap B$. Then the following hold.
\begin{enumerate}
\item[(i)] We have $\Phi\leq\norm{\aut
T}A\cap\norm{\aut T}B$. 
\item[(ii)] The groups $(\norm{\aut T}A\cap \Theta)T/T$ and
$(\norm{\aut T}B\cap\Theta)T/T$ are conjugate to the subgroup in
column {\rm X} of~\cite[Results Matrix]{Kleidman}, so that 
$$
\pi(\norm{\aut
T}A\cap\Theta)\cong\pi(\norm{\aut T}B\cap\Theta)\cong\Z_2\times\Z_2.
$$
\item[(iii)] We have $\Phi\leq \norm{\aut T}C$ and 
$(\norm{\aut T}C\cap\Theta)T/T$ is conjugate to the subgroup
in column {\rm VII} of~\cite[Results Matrix]{Kleidman}, so that $\pi(\norm{\aut
T}C\cap\Theta)\cong\sy 3$ and
\begin{equation}\label{3eq}
\left|\norm{\aut T}{C}:\norm{\aut T}{\{A,
B\}}\right|=3.
\end{equation}
\item[(iv)] We have $T\norm{\aut
T}{\{A,B\}}=T\Phi\left<\vartheta\right>$, where $\vartheta$ is as in
Lemma~$\ref{factlem}$(ii), so that $\pi(\norm{\aut
T}{\{A,B\}}\cap\Theta)/T\cong\Z_2$. 
\end{enumerate}
\end{lemma}
\begin{proof}
Claims (i)--(ii) 
can easily be verified by inspection of~\cite[Results Matrix]{Kleidman}. 
In~(iii) we only need to prove~\eqref{3eq}. Let $N_1=\norm{\aut T}{C}$ and $N_2=\norm{\aut T}{\{A,
B\}}$. Clearly $N_2\leq N_1$.
By~\cite[Proposition~3.1.1(vi)]{Kleidman}, $C\cong{\sf G}_2(q)$,
and~\cite[Proposition~3.1.1(iii)]{Kleidman} shows that 
$\pi(N_1\cap\Theta)\cong\sy 3$. From~\cite[Results Matrix]{Kleidman} we obtain $\pi(N_2\cap\Theta)=\Z_2$. As in the proof of Lemma~\ref{factlem},
we have $N_1\cap T=N_2\cap T=C$. 
As $T=\ker \pi$ this implies $N_1\cap\ker\pi=N_2\cap\ker\pi$, and so
$|N_1\cap\Theta|=3\cdot|N_2\cap\Theta|$. Since $\Phi\leq N_1\cap N_2$ we
have $N_1\Theta=N_2\Theta=\aut T$, and so $|N_1|=3\cdot|N_2|$, as
required. In~(iv) we notice that $T\Phi\left<\vartheta\right>\leq T\norm{\aut
T}{\{A,B\}}$. On the other hand, (iii)~implies that $|T\Phi\left<\vartheta\right>|=|T\norm{\aut
T}{\{A,B\}}|$, hence equality follows.
\end{proof}

\section{Innately transitive groups with a non-abelian, simple plinth}\label{simsec}

In this section we prove our second main theorem, namely
Theorem~\ref{main}, which is a consequence of the
following result.

\begin{theorem}\label{simth}
Let $G$ be an innately transitive permutation group  on $\Omega$
with a non-abelian, simple
plinth $T$, let $\omega\in\Omega$, $\E\in\cd G$, and let $W$ be the stabiliser of $\E$ in
$\sym\Omega$. Then $|\E|\leq 3$ and the following hold.
\begin{enumerate}
\item[(i)] Suppose that $\E$ is homogeneous. Then $|\E|=2$, 
$W$ is a maximal
subgroup of $\sym\Omega$ or $\Alt\Omega$, and $G$, $T$,  $W$, the
subgroups $K\in\K_\omega(\E)$, and $|\Omega|$ are as in Table~$\ref{table}$.
 In particular, the set $\K_\omega(\E)$ contains two isomorphic
subgroups. Moreover, the group $G$ is quasiprimitive and 
$T$ is the unique minimal normal subgroup of $G$. Moreover exactly one
of the following holds:
\begin{enumerate}
\item[(a)] $|\cd G|=1$;
\item[(b)] $|\cd G|=3$, $T$ is as in row~$3$ of Table~$\ref{table}$, $G\leq
T\Phi$ where $\Phi$ is the group of
field automorphisms of $T$.
\end{enumerate}
\item[(ii)] Suppose that $|\E|=3$.  If 
$W$ is the stabiliser in $\sym\Omega$ of $\E$, then 
$T$,  $W$, the elements of $\K_\omega(\E)$, and $|\Omega|$ are as in Table~$\ref{tablemult}$.
\end{enumerate}
\end{theorem}
\tiny
\begin{center}
\begin{table}[ht]
$$
\begin{array}{|c|c|c|c|c|c|}
\hline
 & G & T & W & K & |\Omega|\\
\hline\hline
1 &\alt 6\leq G\leq{\sf P}\Gamma{\sf L}_2(9)& \alt 6& \sy 6\wr\sy 2& \alt 5 & 36\\
\hline
2 & \mat{12}\leq G\leq\aut{\mat{12}}&\mat{12}& \sy {12}\wr\sy 2 & \mat{11} & 144 \\
\hline
3 & \pomegap 8q\leq G\leq \pomegap 8q\Phi\left<\vartheta\right>
& \pomegap 8q & \sy{(d/2)q^3(q^4-1)}\wr\sy
 2
& \Omega_7(q) & (d^2/4)q^6(q^4-1)^2\\
 &\mbox{$\Phi$: field automorphisms} && d=(4,q^4-1) &&\\
&
\mbox{$\vartheta$ is as in Lemma~\ref{factlem}(ii)} &&&&\\

\hline
4 & \sp 4q\leq G\leq\aut{\sp 4q} & \sp 4q,\ q\geq 4\mbox{ even}
& \sy{q^2(q^2-1)}\wr\sy 2 & \sp 2{q^2}.2 & q^4(q^2-1)^2\\
\hline
\end{array}
$$
\caption{Homogeneous Cartesian decompositions preserved by almost
simple groups}
\label{table}
\end{table}
\end{center}
\begin{center}
 \begin{table} 
$$
 \begin{array}{|c|c|c|c|c|}
 \hline
 & T & W & \K_\omega(\E) & |\Omega| \\
 \hline\hline
 1& \sp {4a}2,\ a\geq 2 &\sy{n_1}\times\sy{n_2}\times\sy{n_3} &
\sp {2a}4\cdot 2,\   {\sf O}^-_{4a}(2),\  {\sf
 O}^+_{4a}(2) & n_1\cdot n_2\cdot n_3\\
& & n_1=|\sp{4a}2:\sp{2a}4\cdot
 2| &&  \\
&& n_2=|\sp{4a}2:{\sf
 O}^-_{4a}(2)|&& \\
&& n_3=|\sp{4a}2:{\sf O}^+_{4a}(2)| && \\
 \hline
 2 & \pomegap 83 &\sy{1080}\times\sy{1120}\times\sy{28431}& \Omega_7(3),\
 \Z_3^6\rtimes \psl 43,\  \pomegap 82 & 34,390,137,600 \\
 \hline
 3 & \sp 62 & \sy{120}\times\sy{28}\times\sy{36}&{\sf G}_2(2),\ {\sf
 O}^-_{6}(2),\ {\sf O}^+_{6}(2) &120,960 \\
  &&\sy{240}\times\sy{28}\times\sy{36}& {\sf G}_2(2)',\ {\sf
 O}^-_{6}(2),\ {\sf O}^+_{6}(2) & 241,920\\
  &&\sy{120}\times\sy{56}\times\sy{36}& {\sf G}_2(2),\ {\sf
 O}^-_{6}(2)',\ {\sf O}^+_{6}(2) & 241,920\\
  &&\sy{120}\times\sy{28}\times\sy{72}& {\sf G}_2(2),\ {\sf
 O}^-_{6}(2),\ {\sf O}^+_{6}(2)' & 241,920\\
 \hline
 \end{array}
 $$
 \caption{Cartesian decompositions with index~3 preserved by almost
 simple groups}\label{tablemult}
 \end{table}
\end{center}
\normalsize
\begin{proof}
Suppose that $\E\in\cd G$. Then Proposition~\ref{minv} implies that 
$T_{(\E)}=T$. Let $\ell$ be the index of $\E$, and 
let $\K_\omega(\E)=\{K_1,\ldots,K_\ell\}$ be the corresponding Cartesian
system for $T$. Then the definition of $\K_\omega(\E)$ implies that
if $\ell\geq 3$ then 
$\{K_1,\ldots,K_\ell\}$ is a strong multiple factorisation of the
finite simple group $T$.
Strong multiple factorisations of
finite simple groups are
defined and classified in~\cite{bad:fact}; in particular it is proved that $\ell\leq
3$.  

(a) If
$\ell=3$ then~\cite[Table~V]{bad:fact} shows that $K_1$, $K_2$, $K_3$
have different sizes. Thus if $\E$ is homogeneous then $\ell=2$ and the factorisation $T=K_1K_2$ is
as in Lemma~\ref{factlem}. Hence $T$, $K_1$, $K_2$, and $|\Omega|$ are as in
Table~\ref{table}. The maximality of $W$ follows from Lemma~\ref{cs2cd}.

Let us now prove that $G$ is quasiprimitive. 
As $T$ is transitive on $\Omega$, we have
$\cent{\sym\Omega}T\cong\norm T{T_\omega}/T_\omega$; see~\cite[Theorem~4.2A]{dm}. On the other hand
$T_\omega=K_1\cap K_2$, and Lemma~\ref{factlem} shows that $\norm T{K_1\cap
K_2}=K_1\cap K_2=T_\omega$. Hence $\cent{\sym\Omega}T=1$, and so $T$ is
the unique minimal normal subgroup of $G$. Hence $G$ is an almost
simple quasiprimitive group acting on $\Omega$.

Now we prove that the information given in the $G$-column of
Table~\ref{table} is correct. Since $T$ is the unique minimal normal
subgroup of $G$, we have that $G$ is an almost simple group and $T\leq
G\leq \aut T$. Let $N=\norm{\aut T}{\{K_1,K_2\}}$. 
Note that $G=TG_\omega$ and
$G_\omega\leq N$. On the other hand, $N$ has the property that, since
$A$ and $B$ are not conjugate in $T$,  
$$
T\cap N=\norm T{K_1}\cap \norm
T{K_2}=K_1\cap K_2=T_\omega, 
$$
and so
the $T$-action on $\Omega$ can be extended to $TN$
with point stabiliser $N$. 
Thus $G\leq TN$. By Lemmas~\ref{factlem}(iv) and~\ref{pomegalem}(iv), for $T\cong
\alt 6$, $\mat{12}$, $\pomegap 8q$, and $\sp 4q$, we have $TN={\sf
P}\Gamma{\sf L}_2(9)$, $\aut{\mat{12}}$, $\pomegap
8q\Phi\left<\vartheta\right>$ (where $\Phi$ is the group of  field
automorphisms and $\vartheta$ is as in Lemma~\ref{factlem}(ii)), 
and $\aut{{\sp 4q}}$, respectively. Hence the assertion follows.

Finally we prove the claim concerning $|\cd G|$.
Suppose that $L_1,\ L_2\leq T$ is such that $|L_1|=|L_2|$, $L_1L_2=T$ and
$L_1\cap L_2=T_\omega$. By~\cite{bad:fact}, the full factorisation $T=K_1K_2$ is unique
up to equivalence, Lemma~\ref{equivlemma}(ii) shows that 
there is an element $\alpha\in\aut T$ such
that $\{K_1,K_2\}^\alpha=\{L_1,L_2\}$, and so $\alpha\in\norm{\aut
T}{T_\omega}=\norm{\aut T}{K_1\cap K_2}$.
Lemma~\ref{factlem}(iii) implies that 
if $T$ is as in row~1, 2, or~4 of Table~\ref{table} then $\norm{\aut T}{\{K_1,K_2\}}=\norm{\aut
T}{T_\omega}$ and so 
$\{L_1,L_2\}=\{K_1,K_2\}^\alpha=\{K_1,K_2\}$. Thus $|\cd G|=1$ in
these cases, as asserted.

Suppose now that  $T\cong\pomegap 8q$ for some $q$. Then we obtain
from  Lemma~\ref{pomegalem}(iii) that
$
\left|\norm{\aut T}{T_\omega}:\norm{\aut T}{\{K_1,K_2\}}\right|=3$,
and so the $\norm{\aut
T}{T_\omega}$-orbit containing $\{K_1,K_2\}$ has 3~elements, which
gives rise to 3~different choices of Cartesian systems with respect to
$\omega$. Let $\E_1$, $\E_2$, and $\E_3$ denote the corresponding
Cartesian decompositions of $\Omega$, such that $\E=\E_1$. We computed
above that $\cent{\sym\Omega}T=1$, and 
this implies
that 
$\norm{\sym\Omega}T=\aut T\cap\sym\Omega$. In other words,
$N=\norm{\sym\Omega}T$ is the largest subgroup of $\aut T$ that extends
the $T$-action on $\Omega$. Since $T$ is a transitive subgroup of $N$,
we have $N=TN_\omega$. As $T_\omega$ is a normal subgroup of
$N_\omega$, it follows that $N\leq T\norm {\aut T}{T_\omega}$. On the
other hand
$$
|T\norm{\aut T}{T_\omega}:\norm{\aut T}{T_\omega}|=|T:T\cap\norm{\aut
T}{T_\omega}|=|T:\norm T{T_\omega}|=|T:T_\omega|,
$$
by Lemma~\ref{factlem}(iii). This shows that the $T$-action on
$\Omega$ can be extended to $T\norm{\aut T}{T_\omega}$ with point
stabiliser $\norm{\aut T}{T_\omega}$. In other words $T\norm{\aut
T}{T_\omega}$ is the largest subgroup of $\aut T$ that extends the
$T$-action on $\Omega$. 
The stabiliser of $\E_1$ in $T\norm{\aut
T}{T_\omega}$ is $T\norm{\aut T}{\{K_1,K_2\}}$. Hence if $G\leq \aut
T$ is such that $T\leq G$ and $G$ leaves the Cartesian decomposition
$\E_1$ invariant, then $G\leq T\norm{\aut
T}{\{K_1,K_2\}}=T\Phi\left<\vartheta\right>$, by Lemma~\ref{pomegalem}(iii). 
If $\cd
G\neq\{\E\}$ then, $G$ leaves $\E_1$, $\E_2$, and $\E_3$
invariant. Therefore $G$ lies in the kernel of the action of
$T\norm{\aut T}{T_\omega}$ on $\{\E_1,\E_2,\E_3\}$. Hence $G\leq
T\Phi$, as required.

(b) Suppose that $|\E|=3$. Then $\{K_1,K_2,K_3\}$ is a strong multiple
factorisation of $T$. Therefore using Table~V in~\cite{bad:fact} we
obtain that $T$, $K_1$, $K_2$, $K_3$, and the degree $|\Omega|=|T:K_1\cap K_2\cap
K_3|$ of $G$ are as in Table~\ref{tablemult}.
\end{proof}

The proof of Theorem~\ref{main} is now easy, because 
Theorem~\ref{simth} implies that
$\cent{\sym\Omega}T=1$, and so $\norm{\sym\Omega}T$ is an almost
simple group with socle $T$. For the proof of
Theorem~\ref{mainver}, notice that $W$ is the full stabiliser of a
Cartesian decomposition $\E$ of $\Omega$. As $G\leq W$, the Cartesian
decomposition $\E$ is also $G$-invariant. Hence Theorem~\ref{simth}
implies the required result.

\section{Acknowledgments}

This paper forms part of an Australian Research Council large
grant project. We are grateful to Cai Heng Li for his valuable
advice. We also wish to thank the anonymous referee for his or her many
suggestions that much improved our exposition: in particular for recommending
that we draw attention to Theorem~\ref{main}, and for suggesting a new version
of Lemma~\ref{cscor}.

\end{document}